# THE RIGHT TIME TO SELL A STOCK WHOSE PRICE IS DRIVEN BY MARKOVIAN NOISE[1]

By Robert C. Dalang and M.-O. Hongler

*Ecole Polytechnique Fédérale de Lausanne*

We consider the problem of finding the optimal time to sell a stock, subject to a fixed sales cost and an exponential discounting rate $\rho$. We assume that the price of the stock fluctuates according to the equation $dY_t = Y_t(\mu\, dt + \sigma \xi(t)\, dt)$, where $(\xi(t))$ is an alternating Markov renewal process with values in $\{\pm 1\}$, with an exponential renewal time. We determine the critical value of $\rho$ under which the value function is finite. We examine the validity of the "principle of smooth fit" and use this to give a complete and essentially explicit solution to the problem, which exhibits a surprisingly rich structure. The corresponding result when the stock price evolves according to the Black and Scholes model is obtained as a limit case.

**1. Introduction.** There are many examples of optimal stopping and optimal control problems in continuous time that involve diffusion processes and that have an explicit solution (see, e.g., [1, 16, 21, 22]), but it is rare that the discrete form of the problem, in which the diffusion is replaced by a random walk, can also be solved explicitly (an exception is [2]; see also [3], Chapter 10). One reason is that in the continuous case, it is possible to use the so-called principle of smooth fit, first studied in detail in [9]; see also [23], [14], Chapters 1 and 6, and [5], as well as [22], page 636, for a discussion of this principle and its history. In the discrete case, the problem is much more combinatorial and no such principle is available.

In this paper, we consider a particular optimal stopping problem in an intermediate situation, in which time is continuous but the driving noise is discrete. We show that the principle of smooth fit holds in some situations but not in others, and that the difference between the situations sheds some

Received November 2002.

[1]Supported in part by a grant from the Swiss National Foundation for Scientific Research.

light on why this principle should hold in the first place. It turns out that this semidiscrete problem admits an essentially explicit solution.

The specific problem we consider is when to sell a stock, subject to a fixed sale cost $a$ and an exponential discounting at rate $\rho$. In the classical Black and Scholes model, the stock price $Y_t$ is a solution of the stochastic differential equation

$$dY_t = Y_t(\mu\, dt + \sigma\, dB_t), \qquad Y_0 = y, \tag{1.1}$$

where $(B_t)$ is a standard Brownian motion, and $\mu$ and $\sigma$ are constants. The problem is to find a stopping time $\tau$ which maximizes the expected reward $E(e^{-\rho\tau}(Y_\tau - a))$. This continuous-time problem can be elegantly solved explicitly (see [18], which was the starting point for this paper).

Here, we consider a semidiscretized form of this problem: The driving noise $dB_t$ is replaced by an alternating renewal process $(\xi(t), t \geq 0)$ with values in $\{\pm 1\}$, with an exponential renewal time with mean $1/\lambda$. Equation (1.1) is replaced by

$$dY_t = Y_t(\mu\, dt + \sigma\xi(t)\, dt), \qquad Y_0 = y,$$

which is a (random) *ordinary* differential equation. The problem is again to find a stopping time $\tau$ which maximizes the expected reward $E(e^{-\rho\tau}(Y_\tau - a))$, and even to find the *value function* $\hat{g}(t, y, s)$, which represents the maximal expected reward if we are at time $t$, the current stock price is $y$, $\xi(t) = s$ and we proceed optimally from time $t$ on.

The Markovian noise process $\xi(t)\, dt$ is sometimes called telegrapher's noise [11] and the process $(Y_t)$ is known as a piecewise deterministic Markov process [4] or a random evolution process [20]. This semidiscrete problem is in a sense "simpler" than the previous one, since it does not appeal to Brownian motion and stochastic differential equations, so the statement of the problem is elementary. The process $\xi(t)$ can be thought of as an up or down trend, which may be appropriate on certain time scales and in certain applications. Furthermore, there is a nontrivial covariance between $\xi(t)$ and $\xi(t+h)$, namely $e^{-2\lambda h}$, which is not the case for white noise $dB_t$. Finally, taking an appropriate limit as $\lambda \uparrow +\infty$ and $\sigma \uparrow +\infty$, we should recover the solution from the Black and Scholes model [although (1.1) will have to be interpreted in the Stratonovitch sense; see Remark 12].

It turns out that the structure of the solution to our problem depends heavily on the relationships between the four parameters $\rho$, $\mu$, $\sigma$ and $\lambda$. If $\rho$ is too small, then it is never optimal to stop and $\sup_\tau E(e^{-\rho\tau}(Y_\tau - a)) = +\infty$. Therefore, the problem is of interest only for sufficiently large values of $\rho$: The critical value is determined in Theorem 1.

It is well known [7] that the optimal stopping rule can be described via a "continuation region" and a "stopping region," and for continuous-time problems, the principle of smooth fit states that the value function $\hat{g}$ should be smooth at the boundary between these two regions. Here, we can establish that the principle of smooth fit does indeed hold at certain such boundary points (see Proposition 7), but not necessarily at others. Therefore, the principle of smooth fit at these points yields necessary conditions on the value function, rather than guesses that need to be confirmed, as is generally the case.

In developing our solution, we essentially had two choices. The first was to "guess" the solution, and then establish that the proposed solution is correct. This is the approach used, for instance, in [16] and [22]. The second choice was to use the general theory of optimal stopping, as developed in [7], to establish existence and basic properties of optimal stopping times. While the first approach is more elementary, it requires heavy algebraic calculations. Therefore, we preferred the second approach: The knowledge that an optimal stopping time exists can be used advantageously to derive conditions that are necessarily satisfied by the value function and this simplifies the calculations (even though they remain intricate).

In the solution to our problem, the possible sign of $\rho - \mu - \sigma$ plays an important role. Intuitively, there are two competing features in the problem: the discounting tends to decrease the effect of the fixed sale cost $a$, which encourages waiting and selling later, versus the discounting of the sale price, which encourages selling immediately. For large values of the sale price, the fixed cost $a$ becomes negligible. When $\rho < \mu + \sigma$, during an up trend, the stock price increases at a faster rate than $\rho$, so no matter what the stock price, it is always optimal not to stop as long as the trend is up. On the other hand, when $\rho > \mu + \sigma$, the discounting is stronger than the increase in stock price, even during an up trend, so for large stock prices, the continuation effect from the fixed sale cost loses out and it becomes optimal to stop, even during an up trend. These observations are confirmed by our analysis.

It turns out that when $\rho - \mu - \sigma$ is negative, we obtain an explicit solution given by algebraic formulas (Theorem 10 for the case where $\mu - \sigma < 0$ and Theorem 14 when $\mu - \sigma > 0$). When $\rho - \mu - \sigma > 0$, the solution is essentially explicit, up to solving one transcendental equation (see Theorem 13 for the case where $\mu - \sigma < 0$ and Theorem 16 if $\mu - \sigma > 0$). Finally, in Remark 12, we show how to recover the solution to the problem where the stock price is given by (1.1) as the limit, when $\lambda \uparrow \infty$, $\sigma \uparrow \infty$ and $\sigma/\sqrt{\lambda} \to \sigma_0$, of the solution to our semidiscrete problem.

**2. Stating the problem.** Consider a two-state continuous-time Markov chain $(\xi(t),\ t \geq 0)$ with state space $S = \{-1, +1\}$, defined on a complete probability space $(\Omega, \mathcal{F}, P)$. Processes such as this are discussed in most textbooks on stochastic processes [10, 20]. We assume that the infinitesimal parameters of this Markov chain are given by the matrix

$$(2.1) \qquad G = \begin{pmatrix} -\lambda & \lambda \\ \lambda & -\lambda \end{pmatrix},$$

where $\lambda > 0$. For $r, s \in S$, let $p_{r,s}(t) = P(\xi(t) = s | \xi(0) = r)$ and set $P(t) = (p_{r,s}(t))$. Then

$$P(t) = \tfrac{1}{2} \begin{pmatrix} 1 + e^{-2\lambda t} & 1 - e^{-2\lambda t} \\ 1 - e^{-2\lambda t} & 1 + e^{-2\lambda t} \end{pmatrix}.$$

Fix positive real numbers $\mu$ and $\sigma$, let

$$V(s) = \mu + s\sigma$$

and consider the process $(Y_t, t \geq 0)$, which is a solution of the equation

$$(2.2) \qquad \frac{dY_t}{dt} = V(\xi(t)) Y_t, \qquad t \geq 0.$$

Equivalently,

$$(2.3) \qquad Y_t = Y_0 \exp\left(\mu t + \sigma \int_0^t \xi(u)\, du\right).$$

Let $(\mathcal{F}_t,\ t \geq 0)$ be the natural filtration of $(\xi(t))$. Clearly, $(\mathcal{F}_t)$ is also the natural filtration of $(Y_t)$. We complete this filtration and then it is also right-continuous. The process $(Y_t)$ is not a Markov process with respect to this filtration, whereas the couple $(Y_t, \xi(t))$ is a Markov process with state space $\mathbb{R}_+ \times S$. We let $P_{y_0, s_0}$ denote the conditional probability given that $Y_0 = y_0$ and $\xi(0) = s_0$, and let $E_{y_0, s_0}$ denote expectation under $P_{y_0, s_0}$.

As mentioned in the Introduction, we assume that $Y_t$ denotes the price of an asset at time $t$ and we wish to sell this asset at the highest possible price, subject to a fixed transaction cost $a > 0$ and a discounting rate $\rho > 0$. That is, the benefit of a sale at time $t$ is given by the reward process $(X_t)$ defined by

$$(2.4) \qquad X_t = e^{-\rho t}(Y_t - a).$$

PROBLEM A.  Find a stopping time $\tau^o$ relative to the filtration $(\mathcal{F}_t)$ such that

$$(2.5) \qquad E_{y_0, s_0}(X_{\tau^o}) = \sup_\tau E_{y_0, s_0}(X_\tau),$$

where the supremum is over all $(\mathcal{F}_t)$-stopping times.

**3. Conditions under which the value function is finite.** The *value function* for Problem A is

$$g(y_0, s_0) = \sup_\tau E_{y_0, s_0}(X_\tau).$$

Of course, Problem A is only interesting if this function is finite. The first theorem identifies the conditions on the parameters of the problem that ensure that this is indeed the case.

THEOREM 1. *Given $y_0 > 0$ and $s_0 \in S$, $g(y_0, s_0) < \infty$ if*

(3.1) $$\rho > \mu - \lambda + \sqrt{\sigma^2 + \lambda^2}.$$

*In fact, condition* (3.1) *holds if and only if $E_{y_0, s_0}(\sup_{t \geq 0} |X_t|) < \infty$, while $\rho < \mu - \lambda + \sqrt{\sigma^2 + \lambda^2}$ if and only if $\sup_{t \geq 0} E_{y_0, s_0}(X_t) = +\infty$.*

PROOF. According to [8], Section III.4, the infinitesimal generator of $(Z_t)$, where $Z_t = (Y_t, \xi(t))$, is the operator $\mathcal{A}$, defined for $f : \mathbb{R}_+ \times S \to \mathbb{R}$ such that $f(\cdot, s)$ is continuously differentiable for each $s \in S$, by

$$\mathcal{A}f(y, s) = V(s) y \frac{\partial f}{\partial y}(y, s) + Gf(y, s),$$

where $Gf(y, s) = G_{s,-1} f(y, -1) + G_{s,+1} f(y, +1)$ and the $G_{s,r}$ are the infinitesimal parameters of $(\xi(t))$ from (2.1).

It is not difficult to check that the law of $Y_t$ under $P_{y_0, s_0}$ is absolutely continuous, with compact support, and we let $p(y_0, s_0; t, y, s)$ denote its density on $\{\xi(t) = s\}$, that is, for all Borel sets $A$,

$$P_{y_0, s_0}\{Y_t \in A, \xi(t) = s\} = \int_A p(y_0, s_0; t, y, s)\, dy.$$

We use this to give a formula for $E_{y_0, s_0}(Y_t)$. Although we could appeal to [13], we prefer, for convenience of the reader, to give the derivation. The Kolmogorov forward equation ([12], Chapter 5.1, page 282) states that for $(y_0, s_0) \in \mathbb{R}_+ \times S$ fixed and all $t > 0$, $(y, s) \in \mathbb{R}_+ \times S$,

(3.2) $$\frac{\partial p}{\partial t}(y_0, s_0; t, y, s) = -\frac{\partial}{\partial y}(V(s) y p(y_0, s_0; t, y, s)) + Gp(y_0, s_0; t, y, s).$$

Set

$$f(y_0, s_0; t, s) = E_{y_0, s_0}(Y_t \mathbb{1}_{\{\xi(t) = s\}}) = \int_0^\infty y p(y_0, s_0; t, y, s)\, dy.$$

Multiply both sides of (3.2) by $y$, then integrate over $[0, \infty[$ with respect to $y$ to find, after an integration by parts and because $y \mapsto p(y_0, s_0; t, y, s)$ has compact support, that

$$\frac{\partial f}{\partial t}(y_0, s_0; t, s) = V(s) f(y_0, s_0; t, s) + Gf(y_0, s_0; t, s).$$

Substitute $s = \pm 1$ and suppress $(y_0, s_0)$ from the notation to get

$$\frac{df}{dt}(t,-1) = (\mu - \sigma)f(t,-1) - \lambda f(t,-1) + \lambda f(t,+1),$$

$$\frac{df}{dt}(t,+1) = (\mu + \sigma)f(t,+1) + \lambda f(t,-1) - \lambda f(t,+1).$$

This is a linear system of two differential equations in the unknowns $t \mapsto f(t, \pm 1)$, governed by a matrix with constant coefficients. The eigenvalues of this matrix are

$$\kappa_1 = \mu - \lambda + \sqrt{\sigma^2 + \lambda^2} \quad \text{and} \quad \kappa_2 = \mu - \lambda - \sqrt{\sigma^2 + \lambda^2}.$$

Therefore,

$$E_{y_0,s_0}(Y_t) = f(y_0, s_0; t, -1) + f(y_0, s_0; t, +1) \sim \exp(\kappa_1 t) \qquad \text{as } t \to \infty$$

and

$$E_{y_0,s_0}(e^{-\rho t} Y_t) \sim \exp((\kappa_1 - \rho)t) \qquad \text{as } t \to \infty.$$

If $\rho < \mu - \lambda + \sqrt{\sigma^2 + \lambda^2}$, then $\kappa_1 - \rho > 0$, so $\sup_{t \geq 0} E_{y_0,s_0}(X_t) = +\infty$ if and only if this inequality holds.

Suppose now that (3.1) holds. To show that $E_{y_0,s_0}(\sup_{t \geq 0} |X_t|) < \infty$, it clearly suffices to show that $E_{y_0,s_0}(\sup_{t \geq 0} M_t) < \infty$, where $M_t = e^{-\rho t} Y_t$. We distinguish two cases.

CASE 1. If $\rho \geq \mu + \sigma$, then $dM_t = (\mu - \rho + \sigma \xi(t)) M_t \leq 0$, so the sample paths of $(M_t)$ are nonincreasing and $\sup_t M_t \leq y_0$, $P_{y_0,s_0}$-a.s. Therefore, $E_{y_0,s_0}(\sup_{t \geq 0} M_t) \leq y_0 < \infty$.

CASE 2. If $\rho < \mu + \sigma$, then we proceed as follows [note first that $\mu - \lambda + \sqrt{\sigma^2 + \lambda^2} < \mu + \sigma$, so the condition $\rho < \mu + \sigma$ is compatible with (3.1)]. The infinitesimal generator $\tilde{\mathcal{A}}$ of the Markov process $Z_t = (M_t, \xi(t))$, which applies to functions $f : \mathbb{R}_+ \times S \to \mathbb{R}$, is given by

$$\tilde{\mathcal{A}} f(y,s) = (\mu - \rho + s\sigma) y \frac{\partial f}{\partial y}(y,s) + G f(y,s).$$

Therefore, the process $(f(M_t, \xi(t)), t \geq 0)$ is a martingale if

$$(3.3) \qquad (\mu - \rho + s\sigma) y \frac{\partial f}{\partial y}(y,s) + G f(y,s) = 0, \qquad (y,s) \in \mathbb{R}_+ \times S.$$

Let $h(z,s) = f(e^z, s)$, so that $h(\cdot, \cdot)$ satisfies the linear equation

$$(\mu - \rho + s\sigma) \frac{\partial h}{\partial z}(z,s) + G h(z,s) = 0, \qquad (z,s) \in \mathbb{R} \times S.$$

Substitute $s = \pm 1$ into this equation to get
$$\frac{\partial}{\partial z}\begin{pmatrix} h(z,-1) \\ h(z,+1) \end{pmatrix} = \begin{pmatrix} \lambda/(\mu-\rho-\sigma) & -\lambda/(\mu-\rho-\sigma) \\ -\lambda/(\mu-\rho+\sigma) & \lambda/(\mu-\rho+\sigma) \end{pmatrix} \cdot \begin{pmatrix} h(z,-1) \\ h(z,+1) \end{pmatrix}.$$
This is a linear system of two differential equations in the unknowns $z \mapsto h(z, \pm 1)$, governed by a matrix with constant coefficients, whose eigenvalues are $0$ and $\tilde{\Omega}$, where
$$\tilde{\Omega} = 2\lambda(\rho-\mu)[\sigma^2 - (\mu-\rho)^2]^{-1}.$$
Note that because $\rho > \mu$ by (3.1) and since we are in Case 2, it follows that $\tilde{\Omega} > 0$. The corresponding eigenvectors are
$$\begin{pmatrix} 1 \\ 1 \end{pmatrix} \quad \text{and} \quad \begin{pmatrix} -(\mu-\rho)+\sigma \\ \mu-\rho+\sigma \end{pmatrix}.$$
Therefore, there are constants $C_1$ and $C_2$ such that
$$h(z,s) = C_1 + C_2(s(\mu-\rho) + \sigma)e^{z\tilde{\Omega}}$$
and so
$$f(y,s) = C_1 + C_2(s(\mu-\rho) + \sigma)y^{\tilde{\Omega}}.$$

Let us choose $C_1 = 0$ and $C_2 = 1$. Given $0 < a < y < b < \infty$, define $T = \inf\{t \geq 0 : M_t \notin [a,b]\}$. The process $(f(M_{t \wedge T}, \xi(t \wedge T)), t \geq 0)$ is a bounded martingale. By (3.1), $\rho > \mu$, so by (2.3) and the law of large numbers,

(3.4) $$\lim_{t \to \infty} M_t = 0, \qquad P_{y,s}\text{-a.s.}$$

and therefore $T < \infty$ $P_{y,s}$-a.s. By the optional stopping theorem ([6], Chapter 4.7, (7.4)), we see that
$$f(y,s) = E_{y,s}[f(M_T, \xi(T))].$$
As $\mu - \rho - \sigma < 0$ and $\mu - \rho + \sigma > 0$, $M_t$ cannot reach $b$ for the first time when $\xi(t) = -1$, or $a$ for the first time when $\xi(t) = +1$, so $\xi(T) = +1$ when $T = b$ and $\xi(T) = -1$ when $T = a$. Let $T_x = \inf\{t \geq 0 : M_t = x\}$. It follows that
$$f(y,s) = f(a,-1)P_{y,s}\{T_a < T_b\} + f(b,+1)P_{y,s}\{T_b < T_a\},$$
which implies that
$$P_{y,s}\{T_b < T_a\} = \frac{f(y,s) - f(a,-1)}{f(b,+1) - f(a,-1)}.$$
Let $a \downarrow 0$ to find that
$$P_{y,s}\left\{\sup_t M_t \geq b\right\} = P_{y,s}\{T_b < \infty\} = \frac{f(y,s)}{f(b,+1)} = \frac{f(y,s)}{\mu-\rho+\sigma}b^{-\tilde{\Omega}}$$
and, therefore, $E_{y,s}(\sup_t M_t) < \infty$ if and only if $\tilde{\Omega} > 1$, that is, $2\lambda(\rho-\mu) > \sigma^2 - (\mu-\rho)^2$, which is equivalent to inequality (3.1) [to see this, isolate the square root in (3.1), square both sides of the inequality and simplify]. Theorem 1 is proved. $\square$

**4. Existence of an optimal stopping time.** In view of Theorem 1, we restrict our study to the situation where condition (3.1) holds, that is, for the remainder of the paper, we make the following assumption.

ASSUMPTION A. The parameters of the problem satisfy

$$\rho > \mu - \lambda + \sqrt{\sigma^2 + \lambda^2}.$$

To apply results from the general theory of optimal stopping in continuous time, we set $X_{+\infty} \equiv 0$.

THEOREM 2. *Under Assumption A, there exists an optimal stopping time $\tau^o$, in other words, $\tau^o$ satisfies* (2.5).

PROOF. We apply Theorem 2.41 of [7]. Notice first that $t \mapsto X_t$ is continuous from $[0, \infty]$ to $\mathbb{R}$. Indeed, the only issue is continuity at $+\infty$, which follows from (3.4) if $\rho < \mu + \sigma$ and from (2.4) and (2.3) if $\rho > \mu + \sigma$, since in this last case,

$$0 \leq e^{-\rho t} Y_t \leq Y_0 \exp((-\rho + \mu + \sigma)t) \qquad \text{for all } t > 0.$$

Furthermore, $(X_t)$ is bounded below by $-a$, adapted and "of class D" [i.e., the family $(X_\tau, \tau$ a stopping time) is uniformly integrable] by Assumption A and Theorem 1. Therefore, the hypotheses of [7], Théorème 2.41, are satisfied and the existence of an optimal stopping time is established. $\square$

**5. First properties of the value function.** From the general theory of optimal stopping in continuous time [7], we know that the solution to Problem A uses *Snell's envelope* of the reward process $(X_t, t \geq 0)$, that is, a supermartingale $(Z_t, t \geq 0)$ such that for all $(y, s)$, $P_{y,s}$-a.s,

$$Z_t = \operatorname{ess\,sup} E_{y,s}(X_\tau | \mathcal{F}_t),$$

where the essential supremum is over all stopping times $\tau \geq t$. Because the reward process has the special form $X_t = e^{-\rho t} f_0(Y_t, \xi(t))$, where

(5.1) $\qquad f_0(y, s) = f_0(y) = y - a \qquad$ (no dependence on $s$),

it follows from [7], Théorème 2.75, that in fact,

$$Z_t = e^{-\rho t} g(Y_t, \xi(t)),$$

where $g(y, s)$ is the value function, and there is an optimal stopping time of the form $\tau^o = \inf\{t \geq 0 : g(Y_t, \xi(t)) = f_0(Y_t)\}$ ([7], Théorème 2.76). We therefore examine properties of the function $g(y, s)$ and of $\{y : g(y, s) = y - a\}$.

PROPOSITION 3. (a) *For $s \in S$, $y \mapsto g(y, s)$ is convex (therefore continuous) and nondecreasing. Furthermore, $g(y, s) \geq \max(f_0(y), 0)$.*

(b) *For $s \in S$, the set $\{y \in \mathbb{R}_+ : g(y, s) = y - a\}$ is an interval $[u_s, +\infty[$ (which may be empty, or, in other words, $u_s = +\infty$ may occur).*

PROOF. (a) We note from (2.3) that the law of $Y_t$ under $P_{y,s}$ is the same as the law of $yY_t$ under $P_{1,s}$ and, therefore,

$$g(y, s) = \sup_\tau E_{1,s}(e^{-\rho\tau}(yY_\tau - a)).$$

For a given stopping time $\tau$ and $s \in S$,

(5.2) $\quad y \mapsto E_{1,s}(e^{-\rho\tau}(yY_\tau - a)) = yE_{1,s}(e^{-\rho\tau}Y_\tau) - aE_{1,s}(e^{-\rho\tau})$

is a nondecreasing and affine function of $y$. Therefore, $y \mapsto g(y, s)$, as the supremum of such functions, is nondecreasing and convex.

Observe that for any $t \geq 0$,

$$g(y, s) \geq E_{y,s}(X_t) = E_{y,s}(e^{-\rho t}Y_t - e^{-\rho t}a) \geq -e^{-\rho t}a.$$

Let $t \to +\infty$ to see that $g(y, s) \geq 0$. It is also clear that $g(y, s) \geq E_{y,s}(X_0) = f_0(y)$, so we conclude that $g(y, s) \geq \max(f_0(y), 0)$.

(b) Let $\mathcal{S}_s = \{y \in \mathbb{R}_+ : g(y, s) = y - a\}$, $s \in \{-1, +1\}$. Then for any $y_0 \in \mathcal{S}_s$ and any stopping time $\tau$,

$$E_{1,s}(e^{-\rho\tau}f_0(y_0 Y_\tau)) \leq g(y_0, s) = y_0 - a.$$

Using (5.1), this is equivalent to

$$y_0(E_{1,s}(e^{-\rho\tau}Y_\tau) - 1) + a(1 - E_{1,s}(e^{-\rho\tau})) \leq 0.$$

The second term is nonnegative, so this inequality implies that $E_{1,s}(e^{-\rho\tau}Y_\tau) - 1 \leq 0$ and, therefore, it remains satisfied for any $y \geq y_0$:

$$y(E_{1,s}(e^{-\rho\tau}Y_\tau) - 1) + a(1 - E_{1,s}(e^{-\rho\tau})) \leq 0.$$

However, this inequality translates back to

$$E_{1,s}(e^{-\rho\tau}f_0(yY_\tau)) \leq y - a.$$

Take the supremum over stopping times $\tau$ to conclude that $g(y, s) \leq y - a$ and so $y \in \mathcal{S}_s$.

We conclude that if $y_0 \in \mathcal{S}_s$ and $y \geq y_0$, then $y \in \mathcal{S}_s$, which shows that either $\mathcal{S}_s = \varnothing$ or $\mathcal{S}_s$ is a semiinfinite interval, as claimed. $\square$

**6. The value function in the continuation region.** For $s \in S$, let $u_s$ be defined as in Proposition 3. We write $u_\pm$ instead of $u_{\pm 1}$. By Proposition 3(b) and [7], Théorèmes 2.18 and 2.45, it is optimal *not to stop* when $(Y_t, \xi(t))$ belongs to the *continuation region*

$$\mathcal{C} = ([0, u_-[ \times \{-1\}) \cup ([0, u_+[ \times \{+1\}),$$

while the smallest optimal stopping time is $\tau^o = \inf\{t \geq 0 : (Y_t, \xi(t)) \in (\mathbb{R}_+ \times S) \setminus \mathcal{C}\}$. Set

$$\hat{g}(t, y, s) = e^{-\rho t} g(y, s).$$

Then the process $(\hat{g}(t, Y_t, \xi(t)), t \geq 0)$ is a supermartingale, while $(\hat{g}(t \wedge \tau^o, Y_{t \wedge \tau^o}, \xi(t \wedge \tau^o)), t \geq 0)$ is a martingale ([7], Théorème 2.75 and (2.12.2)). This has the following consequence.

PROPOSITION 4. (a) *The relationship* $u_- \leq u_+$ *holds.*
(b) *The set* $\{y \in \mathbb{R}_+ : g(y, -1) = y - a\}$ *is nonempty or, in other words,* $u_- < +\infty$.

PROOF. (a) We distinguish two cases, according as $\mu - \sigma < 0$ or not.

CASE 1. $\mu - \sigma < 0$. We show that, in fact, $g(\cdot, -1) \leq g(\cdot, +1)$, which clearly implies $u_- \leq u_+$. Let $\tau_1$ be the first jump time of $(\xi(t))$. Then for any stopping time $\tau$,

$$E_{1,-1}(e^{-\rho\tau} f_0(yY_\tau))$$
$$= E_{1,-1}(e^{-\rho\tau} f_0(yY_\tau) \mathbb{1}_{\{\tau \leq \tau_1\}} + E_{1,-1}(X_\tau | \mathcal{F}_{\tau_1}) \mathbb{1}_{\{\tau_1 < \tau\}})$$
$$\leq E_{1,-1}(e^{-\rho\tau} g(yY_\tau, +1) \mathbb{1}_{\{\tau \leq \tau_1\}} + e^{-\rho\tau_1} g(yY_{\tau_1}, +1) \mathbb{1}_{\{\tau_1 < \tau\}}).$$

Because $\mu - \sigma < 0$, $t \mapsto yY_t$ is nonincreasing on $[0, \tau_1]$ $P_{1,-1}$-a.s. so this is no greater than

$$E_{1,-1}(1 \cdot g(y, +1) \mathbb{1}_{\{\tau \leq \tau_1\}} + 1 \cdot g(y, +1) \mathbb{1}_{\{\tau_1 < \tau\}}) = g(y, +1).$$

Therefore, $g(y, -1) \leq g(y, +1)$.

CASE 2. $\mu - \sigma > 0$. We note that in this case, $t \mapsto Y_t$ is monotone and increasing, and it suffices to consider the case where $u_+ < \infty$. Suppose by contradiction that $u_+ < u_-$. Fix $y \in [u_+, u_-[$, so that $g(y, +1) = y - a$ and, in particular, for any $t \geq 0$,

$$E_{1,+1}(e^{-\rho(\tau_1 \wedge t)}(yY_{\tau_1 \wedge t} - a)) \leq y - a.$$

Let $t_- = (\mu - \sigma)^{-1} \ln(u_-/y)$, so that $t_- > 0$ and on $\{\tau_1 > t_-\}$, $yY_{t_-} = u_-$ $P_{1,-1}$-a.s.

Set $\sigma_1 = \tau_1 \wedge t_-$. Because of the form of the continuation region (given above), $\sigma_1 = \tau^o$, $P_{y,-1}$-a.s. and, therefore,

$$(6.1) \qquad g(y,-1) = E_{1,-1}(e^{-\rho(\tau_1 \wedge t_-)}(yY_{\tau_1 \wedge t_-} - a)).$$

Since

$$Y_{\tau_1 \wedge t_-} = \exp((\mu \pm \sigma)(\tau_1 \wedge t_-)) \qquad P_{1,\pm 1}\text{-a.s.}$$

and the law of $\tau_1$ is exponential with mean $1/\lambda$, both under $P_{1,-1}$ and under $P_{1,+1}$, we see that the right-hand side of (6.1) is bounded above by

$$E_{1,+1}(e^{-\rho(\tau_1 \wedge t_-)}(yY_{\tau_1 \wedge t_-} - a)) \leq g(y,+1) = y - a,$$

because $y \geq u_+$. It follows that $g(y,-1) \leq y - a$, and by Proposition 3(a), this inequality must be an equality. Therefore, $(y,-1)$ belongs to $(\mathbb{R}_+ \times S) \setminus \mathcal{C}$, which contradicts our assumption that $y < u_-$. This proves that $u_- \leq u_+$ as claimed.

(b) If $u_-$ were equal to $+\infty$, then $u_+ = +\infty$ by (a), so the continuation region would be $\mathcal{C} = \mathbb{R}_+ \times S$ and, therefore, $\tau^o \equiv +\infty$ would be the smallest optimal stopping time by [7], Théorème 2.45. However, the reward associated with this stopping time is 0, which is clearly not optimal. $\square$

The supermartingale and martingale properties mentioned just before Proposition 4 translate into

$$(6.2) \qquad \hat{\mathcal{A}}\hat{g}(t,y,s) \leq 0 \qquad \text{for all } (t,y,s) \in \mathbb{R}_+ \times \mathbb{R}_+ \times S$$

and

$$(6.3) \qquad \hat{\mathcal{A}}\hat{g}(t,y,s) = 0 \qquad \text{for } (t,y,s) \in [0,\tau^o[ \times \mathcal{C},$$

where $\hat{\mathcal{A}}$ is the infinitesimal generator of $\hat{Z}_t = (t, Y_t, \xi(t))$. By [8], Section III.4,

$$\hat{\mathcal{A}}\hat{g}(t,y,s) = e^{-\rho t}\left(-\rho g(y,s) + V(s)y\frac{\partial g}{\partial y}(y,s) + Gg(y,s)\right),$$

so

$$(6.4) \qquad -\rho g(y,s) + V(s)y\frac{\partial g}{\partial y}(y,s) + Gg(y,s) = 0, \qquad (y,s) \in \mathcal{C}.$$

A priori, this equality is only satisfied a.e., but since $y \mapsto g(y,s)$ is convex, hence continuous, the fact that (6.4) is satisfied a.e. implies that it is, in fact, satisfied everywhere in $\mathcal{C}$ [let $y_n \to y$ with (6.4) satisfied at each $(y_n, s) \in \mathcal{C}$ and take $y_n$ monotone increasing, then decreasing, to conclude that both the left and right derivative of $g(\cdot, s)$ agree at $(y,s)$].

Equality (6.4) provides us with a specific form for $g|_\mathcal{C}$.

PROPOSITION 5. (a) *Define*

$$\Omega_\pm = \frac{(\lambda + \rho)\mu \pm \sqrt{\lambda^2\mu^2 + \sigma^2(\rho^2 + 2\lambda\rho)}}{\mu^2 - \sigma^2} \tag{6.5}$$

*and*

$$w_\pm = 1 + \frac{\rho}{\lambda} - \frac{\mu + \sigma}{\lambda}\Omega_\pm. \tag{6.6}$$

*There are constants $C_-$ and $C_+$ such that for $0 \leq y \leq u_-$,*

$$g(y, -1) = C_- w_- y^{\Omega_-} + C_+ w_+ y^{\Omega_+}, \tag{6.7}$$

$$g(y, +1) = C_- y^{\Omega_-} + C_+ y^{\Omega_+}. \tag{6.8}$$

(b) *Let*

$$b = \lambda(\lambda + \rho - \mu - \sigma)^{-1} \quad and \quad \Omega = (\lambda + \rho)(\mu + \sigma)^{-1}. \tag{6.9}$$

*There is a constant $C$ such that, for $u_- \leq y \leq u_+$,*

$$g(y, +1) = by - a\frac{\lambda}{\lambda + \rho} + Cy^\Omega. \tag{6.10}$$

PROOF. (a) Equation (6.4), written for $(y, s) \in [0, u_-[ \times S$, gives

$$y(\mu - \sigma)\frac{\partial g}{\partial y}(y, -1) - (\lambda + \rho)g(y, -1) + \lambda g(y, +1) = 0, \tag{6.11}$$

$$y(\mu + \sigma)\frac{\partial g}{\partial y}(y, +1) + \lambda g(y, -1) - (\lambda + \rho)g(y, +1) = 0. \tag{6.12}$$

The same change of variables that was used to solve (3.3) transforms these two equations into a linear system of differential equations governed by a matrix with constant coefficients, whose characteristic polynomial is $(\mu^2 - \sigma^2)^{-1}$ multiplied by

$$p(w) = (\mu^2 - \sigma^2)w^2 - 2\mu(\lambda + \rho)w + (\rho^2 + 2\rho\lambda). \tag{6.13}$$

The roots of this polynomial are easily seen to be $\Omega_-$ and $\Omega_+$ given in (6.5), and the associated eigenvectors are

$$\begin{pmatrix} \omega_- \\ 1 \end{pmatrix} \quad \text{and} \quad \begin{pmatrix} \omega_+ \\ 1 \end{pmatrix},$$

where $\omega_\pm$ are given in (6.6). This leads to the formulas in (6.7) and (6.8).

(b) Equation (6.4), written for $y \in [u_-, u_+]$ and $s = +1$, yields the equation

$$y(\mu + \sigma)\frac{\partial g}{\partial y}(y, +1) - (\lambda + \rho)g(y, +1) + \lambda(y - a) = 0, \tag{6.14}$$

because $g(y, -1) = y - a$ for these $y$. The solution of this first-order linear differential equation is easily seen to be given by (6.10). □

REMARK 6. (a) By isolating the square root in (3.1) and squaring, we see that Assumption A is equivalent to the condition $p(1) > 0$, where $p(\cdot)$ is the polynomial in (6.13).

(b) Assumption A clearly implies $\lambda + \rho - \mu - \sigma > 0$ and, therefore, $\Omega > 1$.

**7. The principle of smooth fit.** Proposition 5 gives the form of the value function in the continuation region, but the numbers of $u_\pm$, $C_\pm$ and $C$ remain to be determined. In many control problems for diffusions [1, 9, 21, 22], this is done using the principle of smooth fit. In the presence of piecewise deterministic processes, it is not a priori clear whether this principle applies, and we will see that this need not be the case. In this problem, the principle of smooth fit states that

$$\lim_{y \uparrow u_\pm} \frac{\partial g}{\partial y}(y, \pm 1) = 1$$

as $\partial f_0 / \partial y \equiv 1$.

PROPOSITION 7. (a) *If $\mu - \sigma > 0$, then the principle of smooth fit is satisfied by $g(\cdot, -1)$ at $u_-$.*

(b) *If $u_+ < +\infty$, then the principle of smooth fit is satisfied by $g(\cdot, +1)$ at $u_+$.*

PROOF. (a) We let $y \uparrow u_-$ in (6.11). Because $g(\cdot, \pm 1)$ is continuous and $g(u_-, -1) = f_0(u_-)$, we find that

$$(7.1) \quad u_-(\mu - \sigma) \lim_{y \uparrow u_-} \frac{\partial g}{\partial y}(y, -1) - (\lambda + \rho) f_0(u_-) + \lambda g(u_-, +1) = 0.$$

For $y \geq u_-$, $\hat{\mathcal{A}}\hat{g}(t, y, -1) \leq 0$ and $\hat{g}(t, y, -1) = e^{-\rho t} f_0(y)$, so

$$y(\mu - \sigma) \frac{df_0}{dy}(y) - (\lambda + \rho) f_0(y) + \lambda g(y, +1) \leq 0.$$

Because $(df_0/dy)(y) \equiv 1$, we let $y \downarrow u_-$ to find that

$$u_-(\mu - \sigma) - (\lambda + \rho) f_0(u_-) + \lambda g(u_-, +1) \leq 0.$$

Since we have assumed that $\mu - \sigma > 0$, together with (7.1), this implies

$$(7.2) \quad \lim_{y \uparrow u_-} \frac{\partial g}{\partial y}(y, -1) \geq 1.$$

By Proposition 3(a), $g(\cdot, -1)$ is convex; therefore, $y \mapsto \frac{\partial g}{\partial y}(y, -1)$ is nondecreasing. As $g(y, -1) = f_0(y)$ for $y \geq u_-$ and $(df_0/dy)(y) \equiv 1$, we conclude that the inequality (7.2) is in fact an equality and (a) is proved.

(b) We let $y \uparrow u_+$ in (6.14) to see, similar to the above, that

$$(7.3) \qquad u_+(\mu + \sigma) \lim_{y \uparrow u_+} \frac{\partial g}{\partial y}(y, +1) - (\lambda + \rho) f_0(u_+) + \lambda f_0(u_+) = 0.$$

For $y \geq u_+$, use the inequality $\hat{\mathcal{A}} \hat{g}(t, y, +1) \leq 0$ to get

$$y(\mu + \sigma) \frac{df_0}{dy}(y) - (\lambda + \rho) f_0(y) + \lambda f_0(y) \leq 0.$$

Since $(df_0/dy)(y) \equiv 1$, we let $y \downarrow u_+$ to find that

$$u_+(\mu + \sigma) - (\lambda + \rho) f_0(u_+) + \lambda f_0(u_+) = 0.$$

Because $\mu + \sigma > 0$, together with (7.3), this implies

$$\lim_{y \uparrow u_+} \frac{\partial g}{\partial y}(y, +1) \geq 1$$

and equality follows as in (a). □

REMARK 8. The formulas in Theorems 10 and 13 below can be used to check that when $\mu - \sigma < 0$, the principle of smooth fit is *not* satisfied by $g(\cdot, -1)$ at $u_-$.

The statements in Proposition 7 can be related to the sample path properties of the process $(Y_t, \xi(t))$. Indeed, when $\mu - \sigma > 0$, $t \mapsto Y_t$ is nondecreasing, so the process $(Y_t, \xi(t))$ can enter the stopping region $[u_-, +\infty[ \times \{-1\}$ through the boundary point $(u_-, -1)$: The principle of smooth fit is satisfied at $(u_-, -1)$ in this case. When $\mu - \sigma < 0$, then $t \mapsto Y_t$ is decreasing on the event $\{\xi(t) = -1\}$, so the only way to enter the region $[u_-, +\infty[ \times \{-1\}$ is if $Y_t > u_-$ and $\xi(t)$ changes from $+1$ to $-1$: In this case, $(Y_t, \xi(t))$ has *not* encountered the boundary point $(u_-, -1)$ and the principle of smooth fit is *not* satisfied there. The same considerations apply at $(u_+, +1)$: The only way for the process $(Y_t, \xi(t))$ to enter the stopping region $[u_+, +\infty[ \times \{+1\}$ is through the boundary point $(u_+, +1)$, and the principle of smooth fit holds at this point.

These observations parallel those in [19], page 304, which Shiryaev pointed out to us shortly before the completion of this paper. A nice feature in the proof of Proposition 7 is that the validity of the principle of smooth fit at the boundary points of the stopping region is a fairly direct consequence of the basic relationships (6.2) and (6.3).

**8. Explicit computation of the value function.** We distinguish four cases, presented below as Theorems 10, 13, 14 and 16, according to the possible relationships between the various parameters. Let $u_\pm$ be as in Proposition 3 and let

$$\Omega_\pm, w_\pm, \Omega, b, C_\pm \quad \text{and} \quad C$$

be as defined in Proposition 5. Only the five numbers $u_\pm$, $C_\pm$ and $C$ remain to be determined, since the other six numbers are given explicitly in Proposition 5. We begin with the following relationships.

LEMMA 9. (a) *The growth of the function* $y \mapsto g(y, +1)$ *as* $y \to \infty$ *is linear.*

(b) *Suppose* $\rho \leq \mu + \sigma$. *Then* $u_+ = +\infty$ *and* $C = 0$.
(c) *Suppose* $\rho > \mu + \sigma$. *Then* $b < 1$, $C > 0$ *and* $u_+ < \infty$.
(d) *Suppose* $\mu - \sigma < 0$. *Then* $\Omega_+ < 0 < \Omega_-$ *and* $C_+ = 0$.
(e) *Suppose* $\mu - \sigma > 0$. *Then* $0 < \Omega_- < \Omega_+$ *and* $\omega_+ < 0 < \omega_-$.

PROOF. (a) In the region $[u_-, +\infty[ \times \{-1\}$, $\hat{\mathcal{A}}\hat{g} \leq 0$ by (6.2) or, equivalently,

$$-\rho(y-a) + y(\mu - \sigma) + \lambda g(y, +1) - \lambda(y-a) \leq 0.$$

This inequality can be written

$$g(y, +1) \leq \frac{\lambda + \rho + \sigma - \mu}{\lambda} y - a \frac{\lambda + \rho}{\lambda}.$$

Therefore, $g(y, +1)$ grows at most linearly when $y \to +\infty$.

(b) Set $\hat{f}(t, y, s) = e^{-\rho t} f_0(y)$ (no dependence on $s$). Observe that for large $y$,

$$\hat{\mathcal{A}}\hat{f}(t, y, +1) = e^{-\rho t}(-\rho f_0(y) + (\mu + \sigma)y)$$
$$= e^{-\rho t}((-\rho + \mu + \sigma)y + a\rho)$$
$$> 0,$$

because we have assumed that $\rho \leq \mu + \sigma$. Assume that $u_+ < \infty$. Then for $y > u_+ \geq u_-$, $\hat{f}(t, y, \pm 1) = \hat{g}(t, y, \pm 1)$, so $\hat{\mathcal{A}}\hat{f}(t, y, +1) = \hat{\mathcal{A}}\hat{g}(t, y, +1) \leq 0$ by (6.2). This contradiction shows that $u_+ = +\infty$.

On the other hand, $\Omega > 1$ by Remark 6(b). Therefore, (6.10) and (a) imply that $C = 0$.

(c) Note that in this case, $C = 0$ is not possible, because $\rho > \mu + \sigma$ implies that $b$ in (6.9) satisfies $b < 1$, so it is not possible to have $by - a\lambda(\lambda + \rho)^{-1} \geq y - a$ for all $y > 0$. Therefore, $C > 0$. Because $\Omega > 1$ by Remark 6(b), we conclude from part (a) and (6.10) that $u_+ < \infty$.

(d) Recall that $\Omega_\pm$ are the roots of the polynomial $p(\cdot)$ in (6.13). When $\mu - \sigma < 0$, the product of the roots is negative and $\Omega_+ < \Omega_-$ by (6.5), so $\Omega_+ < 0 < \Omega_-$. From the fact that $g(\cdot, s)$ is continuous and $g(0, s) = 0$, (6.7) and (6.8) imply that $C_+ = 0$.

(e) Because $\mu - \sigma > 0$, the sum and product of the roots of $p(\cdot)$ in (6.13) are positive, and $\Omega_- < \Omega_+$ by (6.5), so $0 < \Omega_- < \Omega_+$. From (6.6), this immediately implies that $\omega_+ - \omega_- < 0$. To get the more precise result in the statement of the lemma, use (6.6) to check that $\omega_+ < 0 < \omega_-$ is equivalent to $\Omega_- < (\lambda + \rho)/(\mu + \sigma) < \Omega_+$, and this follows from the fact that $p((\lambda + \rho)/(\mu + \sigma)) = -\lambda^2$ ($< 0$), as is easily checked. □

THEOREM 10. *Under Assumption* A, *assume that*

$$\rho \leq \mu + \sigma \quad and \quad \mu - \sigma < 0.$$

*Then the value function $g(\cdot, \pm 1)$ (see sketch in Figure* 1*), expressed by the formulas in Proposition* 5*, is characterized by $C_+ = C = 0$, $u_+ = +\infty$ and*

$$(8.1) \qquad C_- = \frac{u_- - a}{\omega_- u_-^{\Omega_-}}, \qquad u_- = a \left[ \frac{1 - (\lambda/(\lambda + \rho))\omega_-}{1 - b\omega_-} \right].$$

PROOF. By Lemma 9(b), $u_+ = +\infty$ and $C = 0$. By Lemma 9(d), $C_+ = 0$. The two remaining unknowns, namely $C_-$ and $u_-$, are determined by matching the value of $g(u_-, -1)$, as expressed in (6.7), with $u_- - a$ and matching $g(u_-, +1)$, as expressed in (6.8), with $g(u_-, +1)$ as expressed in (6.10). This yields, respectively, the relationships

$$C_- \omega_- u_-^{\Omega_-} = u_- - a$$

and

$$(8.2) \qquad C_- u_-^{\Omega_-} = b u_- - a \frac{\lambda}{\lambda + \rho}.$$

Solving the above set of equations for $C_-$ and $u_-$ gives the expressions in the statement of the theorem. □

REMARK 11. The number $\omega_\pm$ can be written

$$(8.3) \qquad \omega_\pm = -\frac{\sigma}{\lambda(\mu - \sigma)} \left[ (\lambda + \rho) \pm \sqrt{(\lambda + \rho)^2 + \lambda^2 \left( \frac{\mu^2}{\sigma^2} - 1 \right)} \right]$$

and this clearly implies that $w_- > 0$, because $\mu - \sigma < 0$. Therefore, (8.1) implies $u_- > a$ and $C_- > 0$, as expected.

REMARK 12 (The white noise limit). Fix $\sigma_0 > 0$ and set $\sigma = \sigma_0 \sqrt{\lambda}$, so that the process $(Y_t)$ satisfies

$$(8.4) \qquad dY_t = Y_t(\mu \, dt + \sigma_0 \sqrt{\lambda} \xi(t) \, dt), \qquad Y_0 = y.$$

Observe that the covariance of $\sqrt{\lambda}\xi(t)$ and $\sqrt{\lambda}\xi(t+h)$ is $\lambda e^{-2\lambda h}$, and we easily check that the noise source $\sqrt{\lambda}\xi(t)$ converges to Gaussian white noise when $\lambda \to +\infty$. In addition, the solution $(Y_t)$ of (8.4) converges weakly [24] to the diffusion process that satisfies

$$(8.5) \qquad dZ_t = Z_t(\mu\, dt + \sigma_0\, dW_t), \qquad Z_0 = y_0,$$

where the stochastic differential equation (SDE) has to be interpreted in the Stratonovich sense. When the stock price is governed by this diffusion equation, the solution to the optimal stopping problem is well known. If we rewrite the SDE (8.5) as the Itô SDE

$$dZ_t = Z_t((\mu + \sigma_0^2/2)\, dt + \sigma_0\, dW_t),$$

then we can use the formulas from [18], Example 10.16, to get the continuation region $[0, u]$ and the value function $g(y)$:

$$
(8.6) \qquad u = a\frac{\Omega_0}{\Omega_0 - 1}, \qquad g(y) = \begin{cases} (u-a)\left[\dfrac{y}{u}\right]^{\Omega_0}, & \text{for } y \leq u, \\ u - a, & \text{for } y \geq u, \end{cases}
$$

with

$$
\Omega_0 = \frac{-\mu + \sqrt{\mu^2 + 2\rho\sigma_0^2}}{\sigma_0^2}.
$$

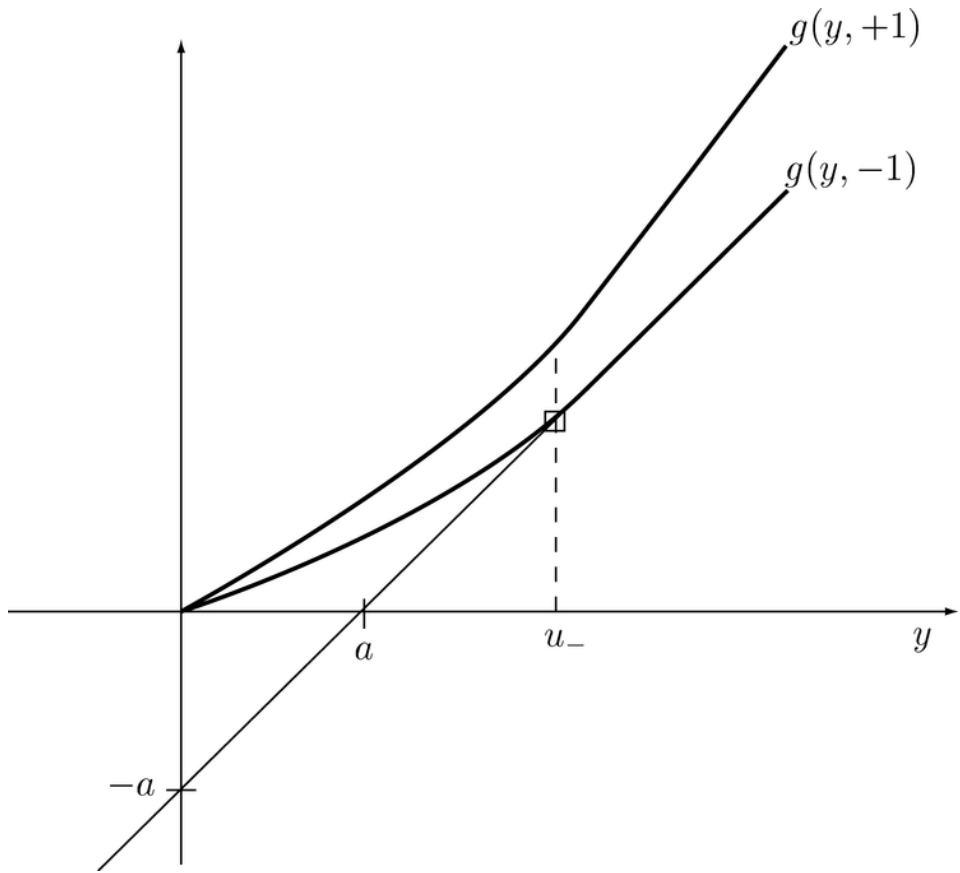

FIG. 1. *The functions $g(\cdot, \pm 1)$ under the hypotheses of Theorem 10. The box indicates the absence of smooth fit.*

The quantities $\Omega_\pm(\lambda)$ and $\omega_\pm(\lambda)$ related to (8.4) are obtained by replacing $\sigma$ in (6.5) and (6.6) with $\sigma_0\sqrt{\lambda}$, which gives

$$
\Omega_\pm(\lambda) = \frac{(\lambda + \rho)\mu \pm \sqrt{\lambda^2\mu^2 + \lambda\sigma_0^2(\rho^2 + 2\lambda\rho)}}{\mu^2 - \lambda\sigma_0^2},
$$

$$
\omega_\pm(\lambda) = 1 + \frac{\rho}{\lambda} - \frac{\mu + \sigma_0\sqrt{\lambda}}{\lambda}\Omega_\pm(\lambda).
$$

Therefore,

$$
\lim_{\lambda \to +\infty} \Omega_-(\lambda) = \Omega_0 \qquad \text{and} \qquad \lim_{\lambda \to +\infty} \omega_\pm(\lambda) = 1.
$$

Furthermore, $u_-$ in (8.1) can be written in the form

$$
u_-(\lambda) = a\frac{\lambda + \rho - \mu - \sigma_0\sqrt{\lambda}}{\lambda + \rho}\left[\frac{\Omega_-(\lambda)}{\Omega_-(\lambda) - 1}\right],
$$

which converges, as $\lambda \to +\infty$, to $u$ as given in (8.6). The general theory of [15] and [17] predicts that $\lim_{\lambda \to +\infty} g(y, \pm 1) = g(y)$, but we can also easily check this from (8.1) and the formulas of Proposition 5.

THEOREM 13. *Assume that*

$$\rho > \mu + \sigma \quad \text{and} \quad \mu - \sigma < 0.$$

*(Note that Assumption* A *is necessarily satisfied in this case.) Then the value function $g(\cdot, \pm 1)$ (see the sketch in Figure 2), expressed by the formulas in Proposition 5, is characterized by*

$$
(8.7) \qquad C_+ = 0 \quad and \quad C_- = \frac{u_- - a}{\omega_- u_-^{\Omega_-}},
$$

*where $u_-$ is the smallest solution of the transcendent equation*

$$
u_- - \frac{C\omega_-}{1 - b\omega_-}u_-^\Omega = \hat{a},
$$

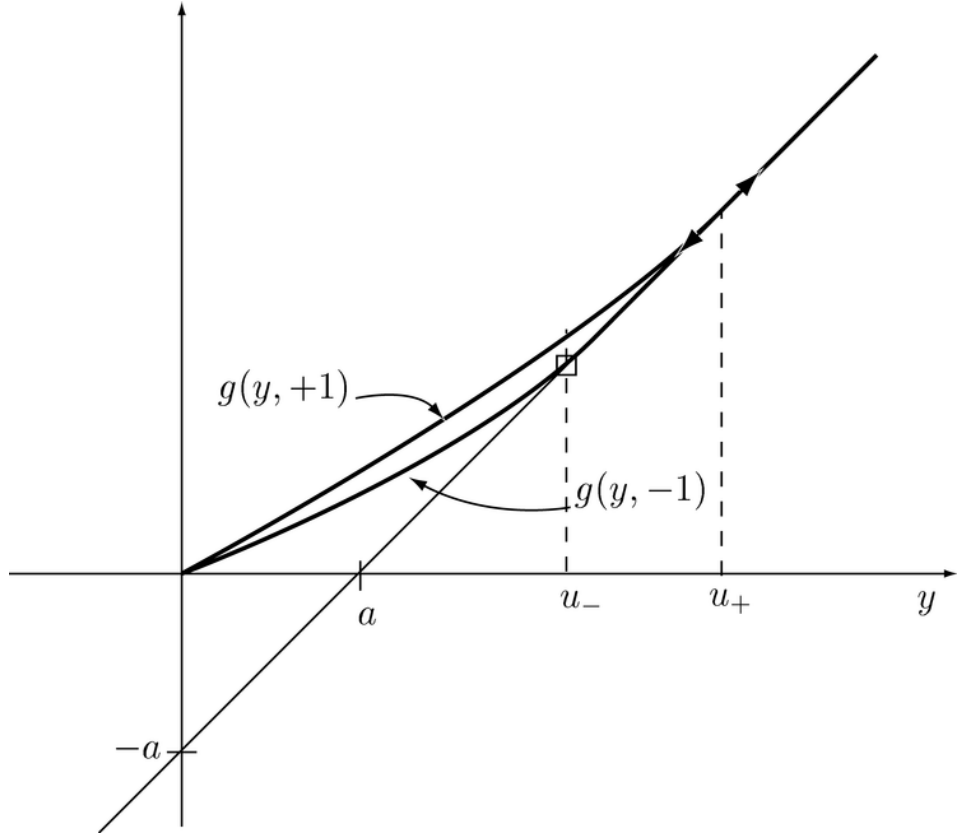

FIG. 2. *The functions $g(\cdot, \pm 1)$ under the hypotheses of Theorem 13. The box (resp. arrow) indicates the absence (resp. presence) of smooth fit.*

$b$ and $\Omega$ are given in (6.9), $\hat{a}$ denotes the expression on the second right-hand side of (8.1),

$$\text{(8.8)} \qquad C = \frac{1-b}{\Omega u_+^{\Omega-1}} \quad \text{and} \quad u_+ = a \frac{\rho}{\rho - \mu - \sigma}.$$

PROOF. By Lemma 9(d), $C_+ = 0$. By Lemma 9(c), $u_+ < \infty$. The continuity of $g(\cdot, -1)$ at $u_-$ and of $g(\cdot, +1)$ at $u_-$ and $u_+$ give, respectively,

$$\text{(8.9)} \qquad C_- \omega_- u_-^{\Omega_-} = u_- - a,$$

$$\text{(8.10)} \qquad C_- u_-^{\Omega_-} = bu_- - a\frac{\lambda}{\lambda + \rho} + C u_-^{\Omega}$$

and

$$bu_+ - a\frac{\lambda}{\lambda+\rho} + Cu_+^\Omega = u_+ - a. \tag{8.11}$$

By Proposition 7(b), there is a smooth fit of $g(\cdot,+1)$ at $u_+$. Accordingly,

$$b + C\Omega u_+^{\Omega-1} = 1. \tag{8.12}$$

The last equation furnishes the formula for $C$, which, when plugged into (8.11), gives after some simplification the formula for $u_+$.

From (8.9), we obtain directly the expression for $C_-$ given in (8.7). Plugging the formulas for $C_-$ and $C$ into (8.10) yields

$$\Psi(u_-) = \hat{a} > a, \qquad \text{where } \Psi(u) = u - \frac{C\omega_-}{1-b\omega_-}u^\Omega, \tag{8.13}$$

with $\hat{a}$ as in the statement of the theorem.

We check that the equation $\psi(u) = \hat{a}$ has, in fact, two solutions and that $u_-$ is the smaller of the two. Observe that

$$\Psi''(u) = -\frac{C\omega_-\Omega(\Omega-1)}{1-b\omega_-}u^{\Omega-2}. \tag{8.14}$$

From Remark 6(b), $\Omega > 1$, and $u_+ > 0$ by (8.8) and the hypothesis $\rho > \mu+\sigma$. Since $b < 1$ by Lemma 9(c), we observe from (8.8) that $C > 0$. By (8.3), $\omega_- > 0$. Therefore, the numerator in (8.14) is positive. A direct calculation using the formula for $\omega_-$ in (6.6) shows that

$$1 - b\omega_- = \frac{\mu+\sigma}{\lambda}b(\Omega_- - 1).$$

Because $p(1) > 0$ by Remark 6(a) and $\Omega_+ < 0 < \Omega_-$ by Lemma 9(d), it follows that, in fact, $\Omega_- > 1$ and, therefore, the denominator in (8.14) is positive.

The above shows that $\Psi''(u) < 0$, so $\Psi$ is strictly concave. From (8.13) and the fact that $\Omega > 1$, it follows that $\Psi(0) = 0$, $\Psi(u) < u$ for all $u > 0$ and $\lim_{u\to+\infty}\Psi(u) = -\infty$. This implies that the equation $\Psi(u) = \hat{a}$ has zero, one or two solutions. Since $u_-$ is a solution, one of the last two occurs. No solution can be less than $a$, since $\hat{a} > a$.

To check that there are exactly two solutions of the equation $\Psi(u) = \hat{a}$ and that $u_-$ is the smaller of the two, we observe from (8.13) that

$$\Psi(u_+) = u_+ - \frac{\omega_-}{1-b\omega_-}Cu_+^\Omega. \tag{8.15}$$

By (8.11),

$$Cu_+^\Omega = (1-b)u_+ - a\frac{\rho}{\lambda+\rho}.$$

Replace $Cu_+^\Omega$ in (8.15) by the right-hand side above, to find, after simplification, that

$$\Psi(u_+) = \frac{1}{1 - b\omega_-}\left(u_+(1 - \omega_-) + \omega_- a \frac{\rho}{\lambda + \rho}\right).$$

Therefore, $\Psi(u_+) > \hat{a}$, since this inequality is now equivalent to $u_+(1-\omega_-) > a(1 - \omega_-)$, which holds since $u_+ > a$ by (8.8), and $\omega_- < 1$ as we now show. Indeed, by (6.6), $\omega_- < 1$ is equivalent to $\rho - (\mu + \sigma)\Omega_- < 0$, which, by (6.5), is in turn equivalent to

$$\rho < \frac{1}{\mu - \sigma}(\mu(\lambda + \rho) - \sqrt{\lambda^2\mu^2 + \sigma^2(\rho^2 + 2\lambda\rho)}\,).$$

Multiply both sides by $\mu - \sigma$, changing the direction of the inequality since $\mu - \sigma < 0$, isolate the square root, then square both sides and simplify to see that this inequality reduces to $\sigma^2 > \mu^2$, which is satisfied by hypothesis.

The inequality $\Psi(u_+) > \hat{a}$ and the properties of $\Psi$ mentioned above imply that one of the solutions of the equation $\Psi(u) = \hat{a}$ is larger than $u_+$ and the other, which is smaller than $u_+$ is therefore $u_-$. □

THEOREM 14. *Under Assumption* A *assume that*

$$\rho \leq \mu + \sigma \quad \text{and} \quad \mu - \sigma > 0.$$

*Then the value function* $g(\cdot, \pm 1)$ *(see the sketch in Figure* 3*), expressed by the formulas in Proposition* 5*, is characterized by* $C = 0$, $u_+ = +\infty$,

(8.16) $\quad C_+ = u_-^{-\Omega_+} \dfrac{u_-(\Omega_- - 1) - a\Omega_-}{\omega_+(\Omega_- - \Omega_+)}, \qquad C_- = u_-^{-\Omega_-} \dfrac{u_-(\Omega_+ - 1) - a\Omega_+}{\omega_-(\Omega_+ - \Omega_-)},$

(8.17) $\quad u_- = a\dfrac{N}{D},$

*where*

(8.18) $$N = \omega_+\Omega_+ - \omega_-\Omega_- - \frac{\lambda\omega_+\omega_-}{\lambda + \rho}(\Omega_+ - \Omega_-)$$

*and*

(8.19) $$D = \omega_+(\Omega_+ - 1) - \omega_-(\Omega_- - 1) - \frac{\lambda\omega_+\omega_-}{\lambda + \rho - \mu - \sigma}(\Omega_+ - \Omega_-).$$

PROOF. By Lemma 9(b), $u_+ = +\infty$ and $C = 0$. The continuity of $g(\cdot, -1)$ and $g(\cdot, +1)$ at $u_-$ give, respectively,

(8.20) $$C_-\omega_- u_-^{\Omega_-} + C_+\omega_+ u_-^{\Omega_+} = u_- - a,$$

(8.21) $$C_- u_-^{\Omega_-} + C_+ u_-^{\Omega_+} = bu_- - a\frac{\lambda}{\lambda + \rho},$$

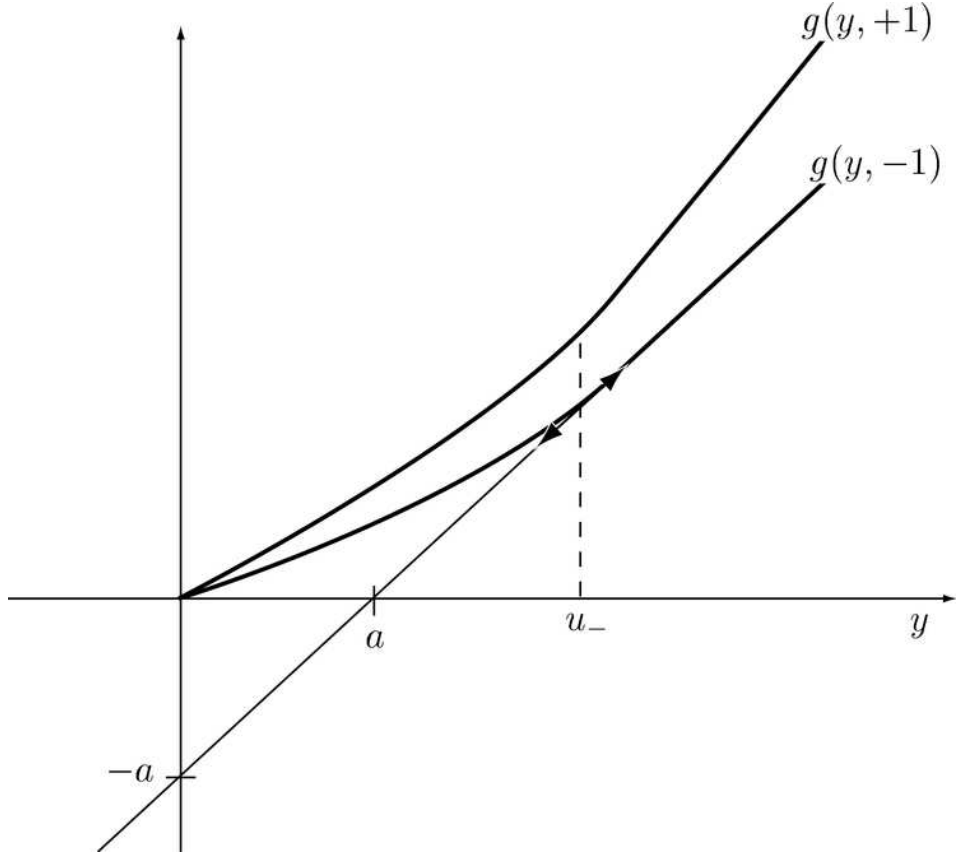

Fig. 3. *The functions $g(\cdot, \pm 1)$ under the hypotheses of Theorem* 14. *The arrow indicates the presence of smooth fit.*

where $b$ is defined in (6.9), and smooth fit of $g(\cdot, -1)$ at $u_-$, which occurs by Proposition 7(a), gives the third equation

$$(8.22) \qquad C_- \omega_- \Omega_- u_-^{\Omega_- - 1} + C_+ \omega_+ \Omega_+ u_-^{\Omega_+ - 1} = 1.$$

Multiply (8.22) by $u_-$, then use (8.20) and (8.22) to express $C_-$ and $C_+$ in terms of $u_-$, which yields the formulas in (8.16). Plug these into (8.21), which becomes a linear equation in $u_-$ and gives (8.17)–(8.19). $\square$

REMARK 15. We note that $u_-$ defined in (8.17) is such that $u_- \geq a$. Indeed, (8.17) can be written

$$u_- = a \frac{N}{D} = a \frac{\chi + \omega_+ - \omega_- + \eta/(\lambda + \rho)}{\chi + \eta/(\lambda + \rho - \mu - \sigma)},$$

with
$$\chi = \omega_+(\Omega_+ - 1) - \omega_-(\Omega_- - 1) \quad \text{and} \quad \eta = \lambda \omega_+ \omega_-(\Omega_- - \Omega_+).$$

We first check that $D < 0$. This is equivalent to verifying
$$\chi < -\eta/(\lambda + \rho - \mu - \sigma).$$

In the definition of $\chi$, of $\eta$ and, on the right-hand side, substitute $\omega_\pm$ from (6.6) to see that this is equivalent to
$$(\lambda + \rho - (\mu + \sigma)\Omega_-)(\Omega_+ - 1) - (\lambda + \rho - (\mu + \sigma)\Omega_-)(\Omega_- - 1)$$
$$< \frac{\lambda^2 \omega_- \omega_+ (\Omega_+ - \Omega_-)}{\lambda + \rho - \mu - \sigma}.$$

With a few algebraic manipulations, $\Omega_+ - \Omega_-$, which is positive by Lemma 9(e), can be factored out on the right-hand side, leading to
$$\lambda + \rho + \mu + \sigma - (\mu + \sigma)(\Omega_+ + \Omega_-) < \frac{\lambda^2 \omega_- \omega_+}{\lambda + \rho - \mu - \sigma}.$$

Now plug into the right-hand side formulas (6.6) for $\omega_\pm$ and simplify to get
$$-1 + \Omega_+ + \Omega_- < \Omega_- \Omega_+.$$

Use the fact that $\Omega_\pm$ are the roots of $p(\cdot)$ in (6.13) to see that this is equivalent to
$$\rho - \mu + \lambda > \sqrt{\lambda^2 + \sigma^2},$$

which holds by Assumption A. Hence, $D < 0$ is established.

We now check that $N < 0$. By Lemma 9(e), $\omega_+ - \omega_- < 0$. Therefore, the facts that $D < 0$ and $(\lambda + \rho - \mu - \sigma)^{-1} \geq (\lambda + \rho)^{-1}$ immediately imply $N < 0$.

It follows that the inequality $u_- \geq a$ is equivalent to $N \leq D$, which becomes
$$\omega_+ - \omega_- \leq \eta[(\lambda + \rho - \mu - \sigma)^{-1} - (\lambda + \rho)^{-1}]$$

(notice that the factor in brackets is $> 0$) and, using (6.6),
$$\frac{\mu + \sigma}{\lambda}(\Omega_- - \Omega_+) \leq \lambda \omega_+ \omega_-(\Omega_- - \Omega_+)[(\lambda + \rho - \mu - \sigma)^{-1} - (\lambda + \rho)^{-1}].$$

Since $\Omega_- - \Omega_+ < 0$ and $\omega_+ \omega_- < 0$ by Lemma 9(e), this inequality does indeed hold, so $u_- > a$ as claimed.

THEOREM 16. *Assume that*
$$\rho > \mu + \sigma \quad \text{and} \quad \mu - \sigma > 0.$$

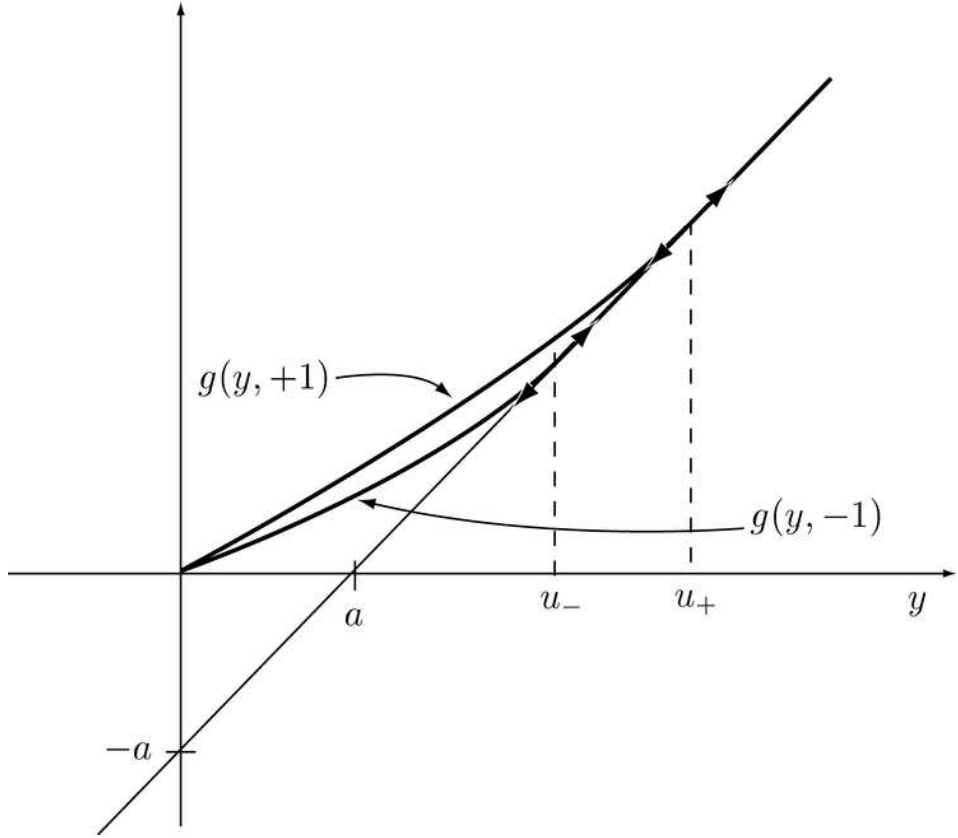

Fig. 4. *The functions $g(\cdot, \pm 1)$ under the hypotheses of Theorem 16. The arrows indicate the presence of smooth fit.*

*(Note that Assumption A is necessarily satisfied in this case.) Then the value function $g(\cdot, \pm 1)$ (see the sketch in Figure 4), expressed by the formulas in Proposition 5, is characterized by*

$$
(8.23) \quad C_+ = u_-^{-\Omega_+} \frac{u_-(\Omega_- - 1) - a\Omega_-}{\omega_+(\Omega_- - \Omega_+)}, \qquad C_- = u_-^{-\Omega_-} \frac{u_-(\Omega_+ - 1) - a\Omega_+}{\omega_-(\Omega_+ - \Omega_-)},
$$

$$
C = \frac{1-b}{\Omega u_+^{\Omega-1}}, \qquad u_+ = a\frac{\rho}{\rho - \mu - \sigma},
$$

*where $b$ and $\Omega$ are given in (6.9) and $u_-$ is the smallest solution of the transcendent equation*

$$
(8.24) \quad \Phi(u) = a\frac{N}{D},
$$

*with $N$ and $D$ defined in* (8.18) *and* (8.19), *and* $\Phi(u)$ *is defined by*

(8.25) $$\Phi(u) \stackrel{def}{=} u - C\frac{\omega_+\omega_-(\Omega_+ - \Omega_-)}{D}u^\Omega.$$

PROOF. By Lemma 9(c), $u_+ < \infty$. The continuity of $g(\cdot, -1)$ and $g(\cdot, +1)$ at $u_-$ give, respectively,

(8.26) $$C_-\omega_- u_-^{\Omega_-} + C_+\omega_+ u_-^{\Omega_+} = u_- - a,$$

(8.27) $$C_- u_-^{\Omega_-} + C_+ u_-^{\Omega_+} = bu_- - a\frac{\lambda}{\lambda + \rho} + Cu_-^\Omega,$$

continuity of $g(\cdot, +1)$ at $u_+$ gives

(8.28) $$bu_+ - a\frac{\lambda}{\lambda + \rho} + Cu_+^\Omega = u_+ - a$$

and Proposition 7 implies two additional equations: one for the smooth fit of $g(\cdot, -1)$ at $u_-$,

(8.29) $$C_-\omega_-\Omega_- u_-^{\Omega_- - 1} + C_+\omega_+\Omega_+ u_-^{\Omega_+ - 1} = 1,$$

and one for the smooth fit of $g(\cdot, +1)$ at $u_+$,

(8.30) $$b + C\Omega u_+^{\Omega - 1} = 1.$$

Observe that equations (8.26) and (8.29) are, respectively, identical to (8.20) and (8.22), which give the formulas for $C_+$ and $C_-$ as in (8.16). Equations (8.28) and (8.30) are, respectively, identical to (8.11) and (8.12), and this gives the formulas for $C$ and $u_+$ as in (8.8). Plug the formulas for $C_\pm$ and $C$ into (8.27) to see that $u_-$ solves the equation

$$\Phi(u) = a\frac{N}{D},$$

where $N$ and $D$ are as in (8.18) and (8.19), and $\Phi(u)$ is as in (8.25).

We show that (8.24) has two solutions, the smaller of which is $u_-$. Since $b < 1$ by Lemma 9(c), we observe from (8.23) that $C > 0$. Since $\Omega_+ - \Omega_- > 0$ and $\omega_+\omega_- < 0$ by Lemma 9(e), $D < 0$ as was observed in Remark 15 and $\Omega > 1$ by Remark 6(b), we see that $\Phi(0) = 0$, $\lim_{u \to +\infty} \Phi(u) = -\infty$ and $\Phi''(u) < 0$ for all $u > 0$, so $\Phi(\cdot)$ is strictly concave. Therefore, (8.24) has zero, one or two solutions. Since $u_-$ is a solution, one of the last two cases occurs.

To show that the equation $\Phi(u) = aN/D$ has exactly two solutions, we proceed as in the last part of the proof of Theorem 13: We show that $\Phi(u_+) > aN/D$, which completes the proof.

From (8.25),

(8.31) $$\Phi(u_+) = u_+ - \frac{\omega_+\omega_-(\Omega_+ - \Omega_-)}{D}Cu_+^\Omega.$$

By (8.28),
$$Cu_+^\Omega = (1-b)u_+ - a\frac{\rho}{\lambda+\rho}.$$

Replace $Cu_+^\Omega$ in (8.31) by the right-hand side above to see that
$$\Phi(u_+) = \frac{u_+}{D}(D - \omega_+\omega_-(\Omega_+ - \Omega_-)(1-b)) + \frac{a}{D}\omega_+\omega_-(\Omega_+ - \Omega_-)\frac{\rho}{\lambda+\rho}.$$

Therefore, the inequality $\Phi(u_+) > aN/D$ is equivalent, after using (8.18) and (8.19) and simplifying, to

(8.32) $\quad (u_+ - a)(\omega_+\Omega_+ - \omega_-\Omega_- - \omega_+\omega_-(\Omega_+ - \Omega_-)) + (\omega_- - \omega_+)u_+ < 0.$

Use (6.6) to see that
$$\omega_+\Omega_+ - \omega_-\Omega_- = (\lambda + \rho - (\mu+\sigma)(\Omega_+ + \Omega_-))\frac{\Omega_+ - \Omega_-}{\lambda}$$

and
$$\omega_- - \omega_+ = \frac{\mu+\sigma}{\lambda}(\Omega_+ - \Omega_-),$$

so (8.32) is equivalent to
$$(u_+ - a)(\lambda + \rho - (\mu+\sigma)(\Omega_+ + \Omega_-) - \lambda\omega_+\omega_-) + (\mu+\sigma)u_+ < 0.$$

Plug in the formula for $u_+$ in (8.23) and use the fact that $\Omega_\pm$ are the roots of $p(w)$ in (6.13) to see that this becomes

(8.33) $\quad (\mu+\sigma)\left(\lambda + \rho - \frac{2\mu(\lambda+\rho)}{\mu-\sigma} - \lambda\omega_+\omega_-\right) + (\mu+\sigma)\rho < 0.$

Using (6.6) and (6.5), we check that
$$\omega_+\omega_- = -\frac{\mu+\sigma}{\mu-\sigma},$$

so (8.33) is equivalent to
$$(\lambda+\rho)(\mu-\sigma) - 2\mu(\lambda+\rho) + \lambda(\mu+\sigma) + \rho(\mu-\sigma) < 0,$$

which reduces to $-2\rho\sigma < 0$. This proves that $\Phi(u_+) > aN/D$, and the properties of $\Phi$ mentioned above imply that one of the solutions of the equation $\Phi(u) = aN/D$ is larger than $u_+$ and the other, which is smaller than $u_+$, is therefore $u_-$. $\square$

INSTITUT DE MATHÉMATIQUES
ECOLE POLYTECHNIQUE FÉDÉRALE
1015 LAUSANNE
SWITZERLAND
E-MAIL: robert.dalang@epfl.ch

INSTITUT DE PRODUCTIQUE ET ROBOTIQUE
ECOLE POLYTECHNIQUE FÉDÉRALE
1015 LAUSANNE
SWITZERLAND
E-MAIL: max.hongler@epfl.ch